\documentclass[letterpaper, 10 pt, conference]{ieeeconf}
\IEEEoverridecommandlockouts
\usepackage[space]{cite}

\usepackage[tbtags]{amsmath}
\usepackage{amsfonts}
\usepackage{amsthm}
\usepackage{mathrsfs}
\usepackage{xspace}
\usepackage{xparse}

\usepackage{tabularx}
\usepackage{booktabs}

\usepackage{algorithm}
\usepackage{algpseudocode}
\usepackage{multicol}

\usepackage{accents}

\usepackage{tikz}
\usetikzlibrary{positioning}

\usepackage{graphicx}
\usepackage{epstopdf}

\usepackage{lipsum}

\usepackage{listings}

\definecolor{Blue}{rgb}{0.0, 0.0, 1.0}
\definecolor{DarkGreen}{rgb}{0.0, 0.2, 0.13}
\definecolor{Purple}{rgb}{0.5, 0.0, 0.5}

\newtheorem{defn}{Definition}

\newtheorem*{BT}{Bochner's Theorem}

\algnewcommand\algorithmicinput{\textbf{Input:}}
\algnewcommand\algorithmicoutput{\textbf{Output:}}
\algnewcommand\Input{\item[\algorithmicinput]}%
\algnewcommand\Output{\item[\algorithmicoutput]}%

\makeatletter
\newcommand\fs@booktabsruled{%
  \def\@fs@cfont{\bfseries\strut}\let\@fs@capt\floatc@ruled
  \def\@fs@pre{\hrule height\heavyrulewidth depth0pt \kern\belowrulesep}%
  \def\@fs@mid{\kern\aboverulesep\hrule height\lightrulewidth\kern\belowrulesep}%
  \def\@fs@post{\kern\aboverulesep\hrule height\heavyrulewidth\relax}%
  \let\@fs@iftopcapt\iftrue
}
\makeatother

\floatstyle{booktabsruled}
\restylefloat{algorithm}

\begin{document}

\title{SReachTools Kernel Module: Data-Driven Stochastic Reachability Using Hilbert Space Embeddings of Distributions}

\author{
Adam~J.~Thorpe,~\IEEEmembership{Student~Member,~IEEE,}
Kendric~R.~Ortiz,~\IEEEmembership{Student~Member,~IEEE,} \\
Meeko~M.~K.~Oishi,~\IEEEmembership{Senior Member,~IEEE}%
\thanks{%
    This material is based upon work supported by the National Science Foundation under NSF Grant Number CNS-1836900.  Any opinions, findings, and conclusions or recommendations expressed in this material are those of the authors and do not necessarily reflect the views of the National Science Foundation.
    The NASA University Leadership initiative (Grant \#80NSSC20M0163) provided funds to assist the authors with their research, but this article solely reflects the opinions and conclusions of its authors and not any NASA entity.
    This research was supported in part by the Laboratory Directed Research and Development program at Sandia National Laboratories, a multimission laboratory managed and operated by National Technology and Engineering Solutions of Sandia, LLC., a wholly owned subsidiary of Honeywell International, Inc., for the U.S. Department of Energy’s National Nuclear Security Administration under contract DE-NA-0003525.  The views expressed in this article do not necessarily represent the views of the U.S. Department of Energy or the United States Government.
}
\thanks{A. Thorpe, K. Ortiz, and M. Oishi are with Electrical
  \& Computer Eng., University of New Mexico, Abq., NM.
  Email: {\tt\{ajthor,kendric,oishi\}@unm.edu}.
}
}



\maketitle

\begin{abstract}
We present algorithms for performing data-driven stochastic reachability as an addition to \texttt{SReachTools}, an open-source stochastic reachability toolbox. Our method leverages a class of machine learning techniques known as kernel embeddings of distributions to approximate the safety probabilities for a wide variety of stochastic reachability problems. By representing the probability distributions of the system state as elements in a reproducing kernel Hilbert space, we can learn the ``best fit'' distribution via a simple regularized least-squares problem, and then compute the stochastic reachability safety probabilities as simple linear operations. This technique admits finite sample bounds and has known convergence in probability. We implement these methods as part of \texttt{SReachTools}, and demonstrate their use on a double integrator system, on a million-dimensional repeated planar quadrotor system, and a cart-pole system with a black-box neural network controller. 
\end{abstract}



\section{Introduction}
\label{section: introduction}

Modern control systems incorporate a wide variety of elements which are resistant to traditional modeling techniques. For instance, systems with autonomous or learning components, human-in-the-loop elements, or poorly characterized stochasticity provide a significant challenge for model-based analysis methods.
In practice, model assumptions can be overly simplistic or fail to capture uncertain behavior, and in some cases are simply wrong. 
As such, these scenarios have brought about a need for algorithms which can provide probabilistic guarantees of safety, even when a comprehensive model of the system is unavailable.
Thus, data-driven techniques for safety analysis are widely applicable for systems with poorly characterized dynamics or uncertainties, and provide an inroad for computing the probability of safety for stochastic systems in a model-free manner.

We present an addition to \texttt{SReachTools} \cite{Vinod_Gleason_Oishi_2019} which enables data-driven stochastic reachability, based on a machine learning technique known as conditional distribution embeddings \cite{song2009hilbert, muandet2017kernel}. 
As a nonparametric technique, kernel distribution embeddings use principles from functional analysis to embed probability distributions as elements in a high-dimensional Hilbert space. 
These techniques have applications to Markov models \cite{grunewalder2012modelling}, partially-observable models \cite{song2010hilbert, nishiyama2012hilbert}, and have recently been used to solve stochastic reachability problems \cite{thorpe2019model, thorpe2019stochastic, thorpe2020datadriven}. 
We incorporate an implementation of the algorithms presented in \cite{thorpe2019model, thorpe2019stochastic, thorpe2020datadriven} into \texttt{SReachTools}, enabling data-driven solutions for stochastic reachability problems that are model-free and distribution agnostic. 

These algorithms have several advantages over comparable model-based techniques. First, the techniques admit finite sample bounds \cite{thorpe2020datadriven} and provide convergence guarantees in the infinite sample case \cite{song2009hilbert}. 
Second, the algorithms largely avoid the curse of dimensionality \cite{bellman2015applied}, which can be a significant computational roadblock when solving dynamic programs. 
Without modification, \cite{thorpe2019model} computes the safety probabilities for a stochastic chain of integrators up to ten thousand dimensions, which is beyond the scope of many existing toolsets. 
However, the techniques are also amenable to several approximative speedup techniques \cite{rahimi2008random, yang2015carte, grunewalder2012conditional}, which have been shown to reduce the computational complexity down to log-linear time. These techniques generally rely upon Fourier approximations of kernel functions and random sampling in the frequency domain, as well as approximations of Gaussian random matrices to alleviate the computational burden. 
In \cite{thorpe2019stochastic}, the authors present an application of one of these techniques, known as random Fourier features \cite{rahimi2008random}, to solve a stochastic reachability problem for a million-dimensional system. 

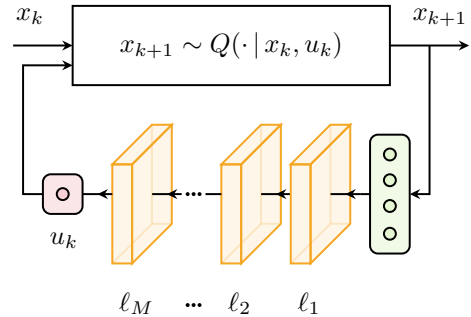
\begin{figure} 
\centering

\tikzset{every picture/.style={line width=0.75pt}} 

\begin{tikzpicture}[x=0.75pt,y=0.75pt,yscale=-1,xscale=1]

\draw  [color={rgb, 255:red, 245; green, 166; blue, 35 }  ,draw opacity=1 ][fill={rgb, 255:red, 245; green, 166; blue, 35 }  ,fill opacity=0.1 ][line width=0.75]  (80,110) -- (90,110) -- (90,160) -- (80,160) -- cycle ;
\draw  [color={rgb, 255:red, 245; green, 166; blue, 35 }  ,draw opacity=1 ][fill={rgb, 255:red, 245; green, 166; blue, 35 }  ,fill opacity=0.1 ] (80,110) -- (100,90) -- (110,90) -- (90,110) -- (80,110) -- cycle ;
\draw  [color={rgb, 255:red, 245; green, 166; blue, 35 }  ,draw opacity=1 ][fill={rgb, 255:red, 245; green, 166; blue, 35 }  ,fill opacity=0.1 ] (110,90) -- (110,140) -- (90,160) -- (90,110) -- (110,90) -- cycle ;

\draw    (115,125) -- (100,125) ;
\draw [shift={(107.5,125)}, rotate = 360] [fill={rgb, 255:red, 0; green, 0; blue, 0 }  ][line width=0.08]  [draw opacity=0] (5.36,-2.57) -- (0,0) -- (5.36,2.57) -- (3.56,0) -- cycle    ;
\draw  [color={rgb, 255:red, 245; green, 166; blue, 35 }  ,draw opacity=1 ][fill={rgb, 255:red, 245; green, 166; blue, 35 }  ,fill opacity=0.1 ][line width=0.75]  (135,110) -- (145,110) -- (145,160) -- (135,160) -- cycle ;
\draw  [color={rgb, 255:red, 245; green, 166; blue, 35 }  ,draw opacity=1 ][fill={rgb, 255:red, 245; green, 166; blue, 35 }  ,fill opacity=0.1 ] (135,110) -- (155,90) -- (165,90) -- (145,110) -- (135,110) -- cycle ;
\draw  [color={rgb, 255:red, 245; green, 166; blue, 35 }  ,draw opacity=1 ][fill={rgb, 255:red, 245; green, 166; blue, 35 }  ,fill opacity=0.1 ] (165,90) -- (165,140) -- (145,160) -- (145,110) -- (165,90) -- cycle ;

\draw    (170,125) -- (155,125) ;
\draw [shift={(162.5,125)}, rotate = 360] [fill={rgb, 255:red, 0; green, 0; blue, 0 }  ][line width=0.08]  [draw opacity=0] (5.36,-2.57) -- (0,0) -- (5.36,2.57) -- (3.56,0) -- cycle    ;
\draw  [color={rgb, 255:red, 245; green, 166; blue, 35 }  ,draw opacity=1 ][fill={rgb, 255:red, 245; green, 166; blue, 35 }  ,fill opacity=0.1 ][line width=0.75]  (170,110) -- (180,110) -- (180,160) -- (170,160) -- cycle ;
\draw  [color={rgb, 255:red, 245; green, 166; blue, 35 }  ,draw opacity=1 ][fill={rgb, 255:red, 245; green, 166; blue, 35 }  ,fill opacity=0.1 ] (170,110) -- (190,90) -- (200,90) -- (180,110) -- (170,110) -- cycle ;
\draw  [color={rgb, 255:red, 245; green, 166; blue, 35 }  ,draw opacity=1 ][fill={rgb, 255:red, 245; green, 166; blue, 35 }  ,fill opacity=0.1 ] (200,90) -- (200,140) -- (180,160) -- (180,110) -- (200,90) -- cycle ;

\draw   (60,30) -- (220,30) -- (220,70) -- (60,70) -- cycle ;
\draw    (220,50) -- (257,50) ;
\draw [shift={(260,50)}, rotate = 180] [fill={rgb, 255:red, 0; green, 0; blue, 0 }  ][line width=0.08]  [draw opacity=0] (5.36,-2.57) -- (0,0) -- (5.36,2.57) -- (3.56,0) -- cycle    ;
\draw    (240,50) -- (240,125) -- (233,125) ;
\draw [shift={(230,125)}, rotate = 360] [fill={rgb, 255:red, 0; green, 0; blue, 0 }  ][line width=0.08]  [draw opacity=0] (5.36,-2.57) -- (0,0) -- (5.36,2.57) -- (3.56,0) -- cycle    ;
\draw    (45,125) -- (35,125) -- (35,60) -- (57,60) ;
\draw [shift={(60,60)}, rotate = 180] [fill={rgb, 255:red, 0; green, 0; blue, 0 }  ][line width=0.08]  [draw opacity=0] (5.36,-2.57) -- (0,0) -- (5.36,2.57) -- (3.56,0) -- cycle    ;
\draw  [fill={rgb, 255:red, 208; green, 2; blue, 27 }  ,fill opacity=0.1 ] (45,119) .. controls (45,116.79) and (46.79,115) .. (49,115) -- (61,115) .. controls (63.21,115) and (65,116.79) .. (65,119) -- (65,131) .. controls (65,133.21) and (63.21,135) .. (61,135) -- (49,135) .. controls (46.79,135) and (45,133.21) .. (45,131) -- cycle ;
\draw  [fill={rgb, 255:red, 208; green, 2; blue, 27 }  ,fill opacity=0.1 ] (52,125) .. controls (52,123.34) and (53.34,122) .. (55,122) .. controls (56.66,122) and (58,123.34) .. (58,125) .. controls (58,126.66) and (56.66,128) .. (55,128) .. controls (53.34,128) and (52,126.66) .. (52,125) -- cycle ;

\draw  [fill={rgb, 255:red, 184; green, 233; blue, 134 }  ,fill opacity=0.2 ] (210,99) .. controls (210,96.79) and (211.79,95) .. (214,95) -- (226,95) .. controls (228.21,95) and (230,96.79) .. (230,99) -- (230,151) .. controls (230,153.21) and (228.21,155) .. (226,155) -- (214,155) .. controls (211.79,155) and (210,153.21) .. (210,151) -- cycle ;
\draw  [fill={rgb, 255:red, 184; green, 233; blue, 134 }  ,fill opacity=0.2 ] (217,105) .. controls (217,103.34) and (218.34,102) .. (220,102) .. controls (221.66,102) and (223,103.34) .. (223,105) .. controls (223,106.66) and (221.66,108) .. (220,108) .. controls (218.34,108) and (217,106.66) .. (217,105) -- cycle ;
\draw  [fill={rgb, 255:red, 184; green, 233; blue, 134 }  ,fill opacity=0.2 ] (217,118) .. controls (217,116.34) and (218.34,115) .. (220,115) .. controls (221.66,115) and (223,116.34) .. (223,118) .. controls (223,119.66) and (221.66,121) .. (220,121) .. controls (218.34,121) and (217,119.66) .. (217,118) -- cycle ;
\draw  [fill={rgb, 255:red, 184; green, 233; blue, 134 }  ,fill opacity=0.2 ] (217,131) .. controls (217,129.34) and (218.34,128) .. (220,128) .. controls (221.66,128) and (223,129.34) .. (223,131) .. controls (223,132.66) and (221.66,134) .. (220,134) .. controls (218.34,134) and (217,132.66) .. (217,131) -- cycle ;
\draw  [fill={rgb, 255:red, 184; green, 233; blue, 134 }  ,fill opacity=0.2 ] (217,145) .. controls (217,143.34) and (218.34,142) .. (220,142) .. controls (221.66,142) and (223,143.34) .. (223,145) .. controls (223,146.66) and (221.66,148) .. (220,148) .. controls (218.34,148) and (217,146.66) .. (217,145) -- cycle ;

\draw    (190,125) -- (205,125) ;
\draw [shift={(197.5,125)}, rotate = 0] [fill={rgb, 255:red, 0; green, 0; blue, 0 }  ][line width=0.08]  [draw opacity=0] (5.36,-2.57) -- (0,0) -- (5.36,2.57) -- (3.56,0) -- cycle    ;
\draw [color={rgb, 255:red, 245; green, 166; blue, 35 }  ,draw opacity=0.2 ]   (100,90) -- (100,140) -- (80,160) ;
\draw [color={rgb, 255:red, 245; green, 166; blue, 35 }  ,draw opacity=0.2 ]   (155,90) -- (155,140) -- (135,160) ;
\draw [color={rgb, 255:red, 245; green, 166; blue, 35 }  ,draw opacity=0.2 ]   (190,90) -- (190,140) -- (170,160) ;
\draw    (80,125) -- (65,125) ;
\draw [shift={(72.5,125)}, rotate = 360] [fill={rgb, 255:red, 0; green, 0; blue, 0 }  ][line width=0.08]  [draw opacity=0] (5.36,-2.57) -- (0,0) -- (5.36,2.57) -- (3.56,0) -- cycle    ;
\draw  [draw opacity=0][fill={rgb, 255:red, 0; green, 0; blue, 0 }  ,fill opacity=1 ] (117,125) .. controls (117,124.45) and (117.45,124) .. (118,124) .. controls (118.55,124) and (119,124.45) .. (119,125) .. controls (119,125.55) and (118.55,126) .. (118,126) .. controls (117.45,126) and (117,125.55) .. (117,125) -- cycle ;
\draw  [draw opacity=0][fill={rgb, 255:red, 0; green, 0; blue, 0 }  ,fill opacity=1 ] (120,125) .. controls (120,124.45) and (120.45,124) .. (121,124) .. controls (121.55,124) and (122,124.45) .. (122,125) .. controls (122,125.55) and (121.55,126) .. (121,126) .. controls (120.45,126) and (120,125.55) .. (120,125) -- cycle ;
\draw  [draw opacity=0][fill={rgb, 255:red, 0; green, 0; blue, 0 }  ,fill opacity=1 ] (123,125) .. controls (123,124.45) and (123.45,124) .. (124,124) .. controls (124.55,124) and (125,124.45) .. (125,125) .. controls (125,125.55) and (124.55,126) .. (124,126) .. controls (123.45,126) and (123,125.55) .. (123,125) -- cycle ;

\draw    (127,125) -- (135,125) ;
\draw    (205,125) -- (210,125) ;
\draw [color={rgb, 255:red, 245; green, 166; blue, 35 }  ,draw opacity=0.2 ]   (100,140) -- (110,140) ;
\draw [color={rgb, 255:red, 245; green, 166; blue, 35 }  ,draw opacity=0.2 ]   (155,140) -- (165,140) ;
\draw [color={rgb, 255:red, 245; green, 166; blue, 35 }  ,draw opacity=0.2 ]   (190,140) -- (200,140) ;
\draw  [draw opacity=0][fill={rgb, 255:red, 0; green, 0; blue, 0 }  ,fill opacity=1 ] (117,181) .. controls (117,180.45) and (117.45,180) .. (118,180) .. controls (118.55,180) and (119,180.45) .. (119,181) .. controls (119,181.55) and (118.55,182) .. (118,182) .. controls (117.45,182) and (117,181.55) .. (117,181) -- cycle ;
\draw  [draw opacity=0][fill={rgb, 255:red, 0; green, 0; blue, 0 }  ,fill opacity=1 ] (120,181) .. controls (120,180.45) and (120.45,180) .. (121,180) .. controls (121.55,180) and (122,180.45) .. (122,181) .. controls (122,181.55) and (121.55,182) .. (121,182) .. controls (120.45,182) and (120,181.55) .. (120,181) -- cycle ;
\draw  [draw opacity=0][fill={rgb, 255:red, 0; green, 0; blue, 0 }  ,fill opacity=1 ] (123,181) .. controls (123,180.45) and (123.45,180) .. (124,180) .. controls (124.55,180) and (125,180.45) .. (125,181) .. controls (125,181.55) and (124.55,182) .. (124,182) .. controls (123.45,182) and (123,181.55) .. (123,181) -- cycle ;

\draw    (30,50) -- (57,50) ;
\draw [shift={(60,50)}, rotate = 180] [fill={rgb, 255:red, 0; green, 0; blue, 0 }  ][line width=0.08]  [draw opacity=0] (5.36,-2.57) -- (0,0) -- (5.36,2.57) -- (3.56,0) -- cycle    ;

\draw (47,143) node [anchor=north west][inner sep=0.75pt]    {$u_{k}$};
\draw (140,50) node    {$x_{k+1} \sim Q( \cdot \,|\, x_{k}, u_{k})$};
\draw (172,173) node [anchor=north west][inner sep=0.75pt]    {$\ell_{1}$};
\draw (137,173) node [anchor=north west][inner sep=0.75pt]    {$\ell_{2}$};
\draw (82,173) node [anchor=north west][inner sep=0.75pt]    {$\ell_{M}$};

\draw (30,30) node [anchor=north west][inner sep=0.75pt]    {$x_{k}$};
\draw (230,30) node [anchor=north west][inner sep=0.75pt]    {$x_{k+1}$};

\end{tikzpicture}
\caption{Feedback control diagram for a stochastic system with a neural network controller.}
\label{fig: neural network diagram}
\end{figure}

Several point-based stochastic reachability techniques are already implemented in \texttt{SReachTools}, based on chance constraints \cite{Vinod_Oishi_2019}, Fourier transforms \cite{Vinod_Oishi_2017}, and particle-based approaches \cite{Lesser_Oishi_Erwin_2013}, among others.
Several existing toolboxes, including Faust${}^{2}$ \cite{soudjani2014faust}, PRISM \cite{kwiatkowska2011prism}, STORM \cite{storm}, and multiple preexisting algorithms in \texttt{SReachTools}, present solutions for stochastic reachability problems and model-checking of continuous and discrete-time Markov chains. 
Unlike most traditional approaches, however, our algorithms are data-driven, meaning we treat the system as a black box, and do not rely upon gridding-based solutions. 
Because of this, we are able to perform stochastic reachability on systems with arbitrary disturbances, as well as systems with autonomous elements such as neural network controllers (Figure \ref{fig: neural network diagram}).
Effectively, this means we can also compute the safety probabilities for autonomous systems and perform neural network verification in a model-free environment. 
Several existing toolboxes, such as Sherlock \cite{dutta2019sherlock}, NNV \cite{nnv}, and Marabou \cite{marabou} tackle the problem of neural network verification, but to the best of our knowledge, our toolbox is the first to be able to compute safety probabilities for stochastic, autonomous systems using backward reachability. 

The paper is organized as as follows:
In Section \ref{section: stochastic reachability}, we present the system model and outline the stochastic reachability problems our algorithms are designed to solve. 
The contribution to \texttt{SReachTools} is presented in Section \ref{section: data-driven stochastic reachability}. 
We present a brief outline of the theory of conditional distribution embeddings and random Fourier features.
Then, we describe the algorithms and their use as part of \texttt{SReachTools}.
In Section \ref{section: numerical experiments}, we present several numerical examples demonstrating the algorithms, including a stochastic integrator system, a repeated planar quadrotor system, and a cart-pole system with a black-box neural network controller. 
Concluding remarks are presented in Section \ref{section: conclusion}.


\section{Stochastic Reachability}
\label{section: stochastic reachability}

We utilize the following notation throughout the paper.
Let $E$ be an arbitrary nonempty space.
The indicator function
$\boldsymbol{1}_{A} : E \rightarrow \lbrace 0, 1 \rbrace$
of $A \subseteq E$ is defined such that
$\boldsymbol{1}_{A}(x) = 1$ if $x \in A$, and
$\boldsymbol{1}_{A}(x) = 0$ if $x \notin A$.
Let $\mathcal{E}$ denote the $\sigma$-algebra on $E$.
If $E$ is a topological space \cite{ccinlar2011probability},
the $\sigma$-algebra generated by the set of all open subsets of $E$
is called the Borel $\sigma$-algebra, denoted by $\mathscr{B}(E)$.
Let $(\Omega, \mathcal{F}, \mathbb{P})$ denote a probability space,
where $\mathcal{F}$ is the $\sigma$-algebra on $\Omega$ and
$\mathbb{P} : \mathcal{F} \rightarrow [0, 1]$ is a \emph{probability measure}
on the measurable space $(\Omega, \mathcal{F})$.
A measurable function $X : \Omega \rightarrow E$ is called a \emph{random variable} taking values in $(E, \mathcal{E})$. The image of $\mathbb{P}$ under $X$, $\mathbb{P}(X^{-1}A)$, $A \in \mathcal{E}$ is called the \emph{distribution} of $X$.

\subsection{System Model}

Consider a Markov control process $\mathcal{H}$ as defined in \cite{Summers_Lygeros_2010}. 
\begin{defn}
    \label{defn: markov control process}
    A Markov control process $\mathcal{H} = (\mathcal{X}, \mathcal{U}, Q)$ is comprised of: 
    \begin{itemize}
        \item 
        $\mathcal{X} \subseteq \mathbb{R}$ a measurable Borel space called the state space;
        
        \item 
        $\mathcal{U} \subset \mathbb{R}$, a compact Borel space called the control space;
        
        \item 
        $Q: \mathscr{B}(\mathcal{X}) \times \mathcal{X} \times \mathcal{U} \rightarrow [0, 1]$, a stochastic kernel that assigns a probability measure $Q(\cdot \,|\, x, u)$ on $(\mathcal{X}, \mathscr{B}(\mathcal{X}))$ to every $(x, u) \in \mathcal{X} \times \mathcal{U}$. 
    \end{itemize}
\end{defn}
The system evolves from an initial condition $x_{0} \in \mathcal{X}$ over a finite time horizon $k = 0, 1, \ldots, N$ with control inputs chosen according to a Markov control policy $\pi$. 
\begin{defn}[Markov Policy]
    A Markov control policy $\pi = \lbrace \pi_{0}, \pi_{1}, \ldots, \pi_{N-1} \rbrace$ is a sequence of universally measurable maps $\pi_{k} : \mathcal{X} \rightarrow \mathcal{U}$, $k = 0, 1, \ldots, N-1$. The set of all admissible Markov policies is denoted as $\mathcal{M}$.
\end{defn}



\subsection{Problem Definitions}

We define the stochastic reachability problems the algorithms can solve as in \cite{Summers_Lygeros_2010} using the stochastic reachability tube definition from \cite{Vinod_Oishi_2020}.

\begin{defn}[Stochastic Reachability Tube]
    \label{defn: stochastic reachability tube}
    Given a finite time horizon $N \in \mathbb{N}$, a stochastic reachability tube is a sequence $\mathcal{A} = \lbrace \mathcal{A}_{0}, \ldots, \mathcal{A}_{N} \rbrace$ of nonempty sets $\mathcal{A}_{k}$, where $k = 0, 1, \ldots, N$.
\end{defn}


\subsubsection{Terminal-Hitting Time Problem}
\label{section: terminal hitting time problem}

Given a target tube $\mathcal{T}$, the terminal-hitting time safety probability $\hat{p}_{x_{0}}^{\pi}$ is defined as the probability that a system following a policy $\pi$ will reach the target set $\mathcal{T}_{N}$ at time $k = N$ while remaining within the target tube $\mathcal{T}$ for all time $k < N$ from an initial condition $x_{0}$.
\begin{equation}
    \hat{p}_{x_{0}}^{\pi}(\mathcal{T}) = \mathbb{E}_{x_{0}}^{\pi} \Biggl[ \Biggl( \prod_{i=0}^{N-1} \boldsymbol{1}_{\mathcal{T}_{i}}(x_{i}) \Biggr) \boldsymbol{1}_{\mathcal{T}_{N}}(x_{N}) \Biggr]
\end{equation}
For a fixed policy $\pi \in \mathcal{M}$, we define the terminal-hitting value functions $W_{k}^{\pi} : \mathcal{X} \rightarrow \mathbb{R}$, $k = 0, 1, \ldots, N$ as 
\begin{align}
    \label{eqn: terminal hitting value functions}
    \begin{split}
        W_{N}^{\pi}(x) &= \boldsymbol{1}_{\mathcal{T}_{N}}(x) \\
        W_{k}^{\pi}(x) &= \boldsymbol{1}_{\mathcal{T}_{k}}(x) 
        \int_{\mathcal{X}} W_{k+1}^{\pi}(y) Q( \mathrm{d} y \,|\, x, \pi_{k}(x))
    \end{split}
\end{align}
Then $W_{0}^{\pi}(x_{0}) = \hat{p}_{x_{0}}^{\pi}(\mathcal{T})$.


\subsubsection{First-Hitting Time Problem}
\label{section: first hitting time problem}

Let $\mathcal{K}$ denote the constraint tube and $\mathcal{T}$ denote the target tube. The first-hitting time safety probability $p_{x_{0}}^{\pi}$ is defined as the probability that a system following a policy $\pi$ will reach the target tube $\mathcal{T}$ at some time $j \leq N$ while remaining within the constraint tube $\mathcal{K}$ for all time $k < j$ from an initial condition $x_{0}$.
\begin{equation}
    p_{x_{0}}^{\pi}(\mathcal{K}, \mathcal{T}) = \mathbb{E}_{x_{0}}^{\pi} \Biggl[ \sum_{j=0}^{N} \Biggl( \prod_{i=0}^{j-1} \boldsymbol{1}_{\mathcal{K}_{i}\backslash\mathcal{T}_{i}}(x_{i}) \Biggr) \boldsymbol{1}_{\mathcal{T}_{N}}(x_{j}) \Biggr]
\end{equation}
For a fixed policy $\pi \in \mathcal{M}$, we define the first-hitting value functions $V_{k}^{\pi} : \mathcal{X} \rightarrow \mathbb{R}$, $k = 0, 1, \ldots, N$ as 
\begin{align}
    \label{eqn: first hitting value functions}
    \begin{split}
        V_{N}^{\pi}(x) &= \boldsymbol{1}_{\mathcal{T}_{N}}(x) \\
        V_{k}^{\pi}(x) &= \boldsymbol{1}_{\mathcal{T}_{k}}(x) 
        + \boldsymbol{1}_{\mathcal{K}_{k}\backslash\mathcal{T}_{k}}(x)
        \int_{\mathcal{X}} V_{k+1}^{\pi}(y) Q( \mathrm{d} y \,|\, x, \pi_{k}(x))
    \end{split}
\end{align}
Then $V_{0}^{\pi}(x_{0}) = p_{x_{0}}^{\pi}(\mathcal{K}, \mathcal{T})$.


\section{Data-Driven Stochastic Reachability}
\label{section: data-driven stochastic reachability}

We consider the case where the stochastic kernel $Q$ is unknown, but observations taken from the system evolution are available. Thus, we have no prior knowledge of the system dynamics or the structure of the disturbance. 
Because $Q$ is unknown, we cannot compute the safety probabilities in \eqref{eqn: first hitting value functions} or \eqref{eqn: terminal hitting value functions} directly.
Instead, we seek to compute an approximation of the safety probabilities by embedding the expectation operator with respect to the stochastic kernel $Q$ in a reproducing kernel Hilbert space and estimating the operator in Hilbert space. 

\begin{defn}[RKHS]
    Let $E$ be an arbitrary, nonempty space and $\mathscr{H}_{E}$ be a Hilbert space over $E$, which is a linear space of functions of the form $f : E \rightarrow \mathbb{R}$. The Hilbert space $\mathscr{H}_{E}$ is a \emph{reproducing kernel Hilbert space} (RKHS) if there exists a positive definite kernel function $k_{E} : E \times E \rightarrow \mathbb{R}$, and it obeys the following properties:
    \begin{subequations}
		\begin{align}
			& k_{E}(x, \cdot) \in \mathscr{H}_{E} && \forall x \in E \\
			\label{eqn: reproducing property}
			& f(x) = \langle f, k_{E}(x, \cdot) \rangle_{\mathscr{H}_{E}} && \forall f \in \mathscr{H}_{E}, \forall x \in E
		\end{align}
	\end{subequations}
	where \eqref{eqn: reproducing property} is called the \emph{reproducing property},
	and for any $x, x' \in E$, we denote $k_{E}(x, \cdot)$ in the RKHS $\mathscr{H}_{E}$ as a function on $E$ such that $x' \mapsto k_{E}(x, x')$.
\end{defn}


We define an RKHS $\mathscr{H}_{\mathcal{X}}$ over $\mathcal{X}$ and $\mathscr{H}_{\mathcal{X} \times \mathcal{U}}$ over $\mathcal{X} \times \mathcal{U}$ with corresponding kernel functions $k_{\mathcal{X}}$ and $k_{\mathcal{X} \times \mathcal{U}}$, and
define the conditional distribution embedding of $Q$ as:
\begin{equation}
    m_{Y|x, u} := \int_{\mathcal{X}} k_{\mathcal{X}}(y, \cdot) Q( \mathrm{d} y \,|\, x, u)
\end{equation}
By the reproducing property, for any function $f \in \mathscr{H}_{\mathcal{X}}$, we can compute the expectation of $f$ with respect to the distribution $Q( \cdot \,|\, x, u)$, where $(x, u) \in \mathcal{X} \times \mathcal{U}$, as an inner product in Hilbert space with the embedding $m_{Y|x, u}$ \cite{thorpe2019model}.
Thus, we can compute the safety probabilities for the first-hitting time problem and the terminal-hitting time problem by computing the expectations in \eqref{eqn: first hitting value functions} and \eqref{eqn: terminal hitting value functions} as Hilbert space inner products.

\subsection{Kernel Distribution Embeddings}

However, we typically do not have access to $m_{Y|x, u}$ directly since the distribution is typically unknown. Instead, we compute an estimate $\hat{m}_{Y|x, u}$ of $m_{Y|x, u}$ using a sample $\mathcal{S} = \lbrace ( x_{i}{}^{\prime}, x_{i}, u_{i} ) \rbrace_{i=1}^{M}$ of $M \in \mathbb{N}$ observations taken i.i.d. from $Q$, where $u_{i} = \pi(x_{i})$ and $x_{i}{}^{\prime} \sim Q(\cdot \,|\, x_{i}, u_{i})$.
The estimate $\hat{m}_{Y|x, u}$ can be found as the solution to a regularized least squares problem:
\begin{equation}
    \label{eqn: regularized least squares}
    \min_{\hat{m}} \frac{1}{M} \sum_{i=1}^{M} \lVert k_{\mathcal{X}}(x_{i}{}^{\prime}, \cdot) - \hat{m}_{Y|x_{i}, u_{i}} \rVert_{\mathscr{H}_{\mathcal{X}}}^{2} + \lambda \lVert \hat{m} \rVert_{\Gamma}^{2}
\end{equation}
where $\lambda > 0$ is the regularization parameter and $\Gamma$ is a vector-valued RKHS.
The solution to \eqref{eqn: regularized least squares} is unique and has the following form:
\begin{equation}
    \label{eqn: estimate}
    \hat{m}_{Y|x, u} = \beta^{\top} \Psi
\end{equation}
where $\beta \in \mathscr{H}_{\mathcal{X}}$ and $\Psi$ is known as a \emph{feature vector}, with elements $\Psi_{i} = k_{\mathcal{X} \times \mathcal{U}}((x_{i}, u_{i}), (x, u))$.
The coefficients $\beta$ are the unique solution to the system of linear equations 
\begin{align}
    \label{eqn: beta}
    (G + \lambda M I) \beta &= \Phi 
\end{align}
where $G = (g_{ij}) \in \mathbb{R}^{M \times M}$ is known as the Gram or kernel matrix, with elements $g_{ij} = k_{\mathcal{X} \times \mathcal{U}}((x_{i}, u_{i}), (x_{j}, u_{j}))$, and $\Phi$ is a feature vector with elements $\Phi_{i} = k_{\mathcal{X}}(x_{i}{}^{\prime}, \cdot)$.
Thus, we can approximate the expectation of a function $f \in \mathscr{H}_{\mathcal{X}}$ using an inner product $\langle \hat{m}_{Y|x, u}, f \rangle_{\mathscr{H}_{\mathcal{X}}}$.
As shown in \cite{thorpe2019model, thorpe2019stochastic}, we can substitute 
the inner product into the backward recursion in \eqref{eqn: first hitting value functions} and \eqref{eqn: terminal hitting value functions} to approximate the safety probabilities.
We have implemented this in \texttt{SReachTools} as \texttt{KernelEmbeddings}. 
We present the algorithm for the first-hitting time problem as Algorithm \ref{algo: kernel embeddings}.


\begin{algorithm}[!t]
    \caption{\texttt{KernelEmbeddings}}
    \label{algo: kernel embeddings}
    
    \begin{algorithmic}[1]
    
        \Input sample $\mathcal{S}$, policy $\pi$, horizon $N$
        \Output value function estimate $V{}_{0}^{\pi}(x)$
        
        \State Initialize
        $V{}_{N}^{\pi}(x) \gets \boldsymbol{1}_{\mathcal{T}_{N}}(x)$
        \State Compute $\hat{m}_{Y|x,\pi(x)}$ 
        \eqref{eqn: estimate}
        using $\mathcal{S}$
        \For{$k \gets N-1$ to $0$}
        \State $
        \mathcal{Y} \gets [V{}_{k+1}^{\pi}(x_{1}{}^{\prime}), \ldots, V{}_{k+1}^{\pi}(x_{M}{}^{\prime})]^{\top}$
        \State
            $
            V{}_{k}^{\pi}(x) \gets
            \boldsymbol{1}_{\mathcal{T}_{k}}(x) +
            \boldsymbol{1}_{\mathcal{K}_{k}\backslash\mathcal{T}_{k}}(x)
            \mathcal{Y}^{\top}
            (G + \lambda M I)^{-1} \Psi
            $
        \EndFor
        \State Return
        $V{}_{0}^{\pi}(x) \approx p_{x}^{\pi}(\mathcal{K}, \mathcal{T})$
        
    \end{algorithmic}

\end{algorithm}

\begin{algorithm}[!t]
    \caption{\texttt{KernelEmbeddingsRFF}}
    \label{algo: kernel embeddings rff}

    \begin{algorithmic}[1]
    
        \Input sample $\mathcal{S}$, policy $\pi$, horizon $N$, sample $\Omega$
        \Output value function estimate $V{}_{0}^{\pi}(x)$
        
        \State Initialize $V{}_{N}^{\pi}(x) \gets \boldsymbol{1}_{\mathcal{T}_{N}}(x)$
        \State Compute RFF approximation 
        of $\hat{m}_{Y|x,\pi(x)}$
        \eqref{eqn: estimate RFF}
        using $\mathcal{S}$ and $\Omega$
        \For{$k \gets N-1$ to $0$}
        \State $\mathcal{Y} \gets [V{}_{k+1}^{\pi}(x_{1}{}^{\prime}), \ldots, V{}_{k+1}^{\pi}(x_{M}{}^{\prime})]^{\top}$
        \State
            $
            V{}_{k}^{\pi}(x) \gets
            \boldsymbol{1}_{\mathcal{T}_{k}}(x) +
            \boldsymbol{1}_{\mathcal{K}_{k}\backslash\mathcal{T}_{k}}(x)
            \mathcal{Y}^{\top} (ZZ^{\top} + \lambda M I)^{-1} Z
            $
        \EndFor
        \State Return $V{}_{0}^{\pi}(x) \approx p_{x}^{\pi}(\mathcal{K}, \mathcal{T})$
        
    \end{algorithmic}
    
\end{algorithm}

\subsection{Random Fourier Features}

Note that in order to compute $\beta$ in \eqref{eqn: beta}, we must solve a system of linear equations that scales with the number of observations $M$ in the sample $\mathcal{S}$.
When $M$ is prohibitively large, 
\cite{rahimi2008random} shows that
by exploiting Bochner's theorem \cite{rudin1962fourier},
we can approximate the kernel function 
by computing its Fourier transform and approximating the Fourier integral in feature space using random realizations of the frequency variable.
\begin{BT} \cite{rudin1962fourier}
  A continuous, translation-invar\-iant kernel function 
  $k(x, x') = \varphi(x - x')$
  is positive definite
  if and only if $\varphi(x - x')$
  is the Fourier transform of a non-negative Borel measure $\Lambda$.
\end{BT}
Following \cite{rahimi2008random}, 
we construct an estimate of the Fourier integral
using a set of $D$ realizations $\Omega = \lbrace \omega_{i} \rbrace_{i=1}^{D}$,
such that $\omega_{i}$ is drawn i.i.d. from the Borel measure $\Lambda$
according to $\omega_{i} \sim \Lambda(\cdot)$.
Then, we can efficiently compute an approximation of the kernel as a sum of cosines.
\begin{equation}
    k(x, x') \approx
    \frac{1}{D}\sum_{i=1}^D \cos(\omega_{i}^{\top}(x - x'))
    \label{eqn: rff integral approximation}
\end{equation}
According to \cite{thorpe2019stochastic}, we can approximate the estimate using \eqref{eqn: rff integral approximation} as: 
\begin{equation}
    \label{eqn: estimate RFF}
    \hat{m}_{Y|x, u} \approx \gamma^{\top} Z 
\end{equation}
where $Z$ is a feature vector computed using the RFF approximation of $k_{\mathcal{X} \times \mathcal{U}}$ and the coefficients $\gamma$ are the solution to the system of linear equations 
\begin{equation}
    \label{eqn: gamma}
    (ZZ^{\top} + \lambda M I) \gamma = \Phi
\end{equation}
This enables a more computationally efficient approximation of the safety probabilities \cite{thorpe2019stochastic} when the number of observations $M$ is large, or when the dimensionality of the system is high. 
We implement this in \texttt{SReachTools} as \texttt{KernelEmbeddingsRFF} and present the algorithm for the first-hitting time problem as Algorithm \ref{algo: kernel embeddings rff}.


\subsection{SReachTools Kernel Module}
\label{section: kernel methods module}

We incorporate the algorithms into \texttt{SReachTools} as modular components, meaning they can be applied to any closed-loop system for which observations are available, and can be applied to either the first-hitting time problem, the terminal-hitting time problem, or the reach-avoid problem \cite{Vinod_Gleason_Oishi_2019}, which can be viewed as a simplification of the terminal-hitting time problem. 

In the implementation, we restrict ourselves to a Gaussian kernel function since the Fourier transform is easy to compute, and to enable a more direct comparison between the two algorithms. The kernel has the form $k(x, x') = \exp(- \lVert x - x' \rVert_{2}^{2}/2 \sigma^{2})$,
where $\sigma$ is known as the bandwidth parameter. 
In order to use the algorithms, we must specify the parameter $\sigma$, the regularization parameter $\lambda$ \eqref{eqn: regularized least squares}, the sample $\mathcal{S}$, and the problem, which is either the first-hitting time problem or the terminal hitting time problem, parameterized by the stochastic reachability tubes $\mathcal{K}$ and $\mathcal{T}$ (Definition \ref{defn: stochastic reachability tube}).

The algorithm parameters $\sigma$ and $\lambda$ are typically chosen via cross-validation, though \cite{Caponnetto_DeVito_2007} suggests a method to select a more optimal rate. The default values are chosen to be $\sigma = 0.1$ and $\lambda = 1$.
We represent a sample $\mathcal{S} = \lbrace (x_{i}{}^{\prime}, x_{i}, u_{i}) \rbrace_{i=1}^{M}$ of size $M \in \mathbb{N}$ taken i.i.d. from a Markov control process $\mathcal{H}$ as a set of matrices, $(X, U, Y)$, where the $i^{\rm th}$ columns of $X$ and $U$ are the observations $x_{i}$ and $u_{i}$, where $u_{i} = \pi(x_{i})$, and the $i^{\rm th}$ column of $Y$ is $x_{i}{}^{\prime}$, where $x_{i}{}^{\prime} \sim Q(\cdot \,|\, x_{i}, u_{i})$.
For \texttt{KernelEmbeddingsRFF}, we are also required to specify the size $D$ of the frequency sample $\Omega$.
The stochastic reachability tubes $\mathcal{K}$ and $\mathcal{T}$ are specified as in \cite{Vinod_Gleason_Oishi_2019}, which is typically defined as a sequence of polyhedra representing the constraints. However, we also include the possibility of representing the tubes via a function-based approach, where the constraints are user-specified indicator functions. 

We then use the sample, the constraint and target tube definitions, and the kernel parameters as inputs to the algorithm and compute the safety probabilities for a point $x_{0} \in \mathcal{X}$ 
using \texttt{SReachPoint} \cite{Vinod_Gleason_Oishi_2019}, 
which returns the safety probability $p_{x_{0}}^{\pi}(\mathcal{K}, \mathcal{T}) \in [0, 1]$.


\section{Numerical Experiments}
\label{section: numerical experiments}

We present several examples to showcase the capabilities of the proposed methods.
Numerical experiments were performed in Matlab on a Intel Xeon CPU with 32 GB RAM, and computation times were obtained using Matlab’s Performance Testing Framework.
Code to reproduce the analysis and figures is available at: 
\texttt{github.com/unm-hscl/ajthor-ortiz-CDC2021}.


\begin{table}[t]
	\caption{Computation times of \texttt{SReachPoint} algorithms for a stochastic double integrator system with $N = 5$.}
	\label{table: computation time}
	\centering
	\normalsize
	\begin{tabular*}{\columnwidth}{
		@{\extracolsep{\fill}} l
		@{\hspace{0.5\tabcolsep}} l
		@{\hspace{0.5\tabcolsep}} c
		}
		\toprule
		
		& 
		& Computation
		\\
		  Algorithm
		& Sample Size
		& Time [s]
		\\
		\midrule
		\texttt{KernelEmbeddings} 
		& $M=2500$
		& $0.72$ s
		\\
		\texttt{KernelEmbeddingsRFF} 
		& $M=2500$, 
		& 
		\\
		& $D=5000$
		& $2.78$ s
		\\
		\midrule
		\texttt{ChanceAffine}
		& 
		& $3.71$ s
		\\
		\texttt{ChanceAffineUniform}
		& 
		& $2.22$ s
		\\
		\texttt{ChanceOpen}
		& 
		& $0.25$ s
		\\
		\texttt{GenzpsOpen}
		& 
		& $0.53$ s
		\\
		\texttt{ParticleOpen}
		& 
		& $0.34$ s
		\\
		\texttt{VoronoiOpen}
		& 
		& $1.69$ s
		\\
		\bottomrule
	\end{tabular*}
\end{table}


\begin{figure}
    \centering
    \includegraphics[keepaspectratio]{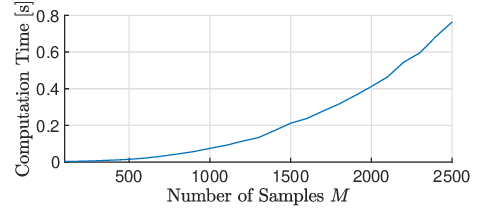}
    \caption{Computation time of \texttt{KernelEmbeddings} as a function of the sample size $M$. The computation time increases roughly exponentially as a the sample size increases.}
    \label{fig: computation time}
\end{figure}


\subsection{Stochastic Chain of Integrators}

\begin{figure*}
    \centering
    \includegraphics[width=\linewidth,height=100pt]{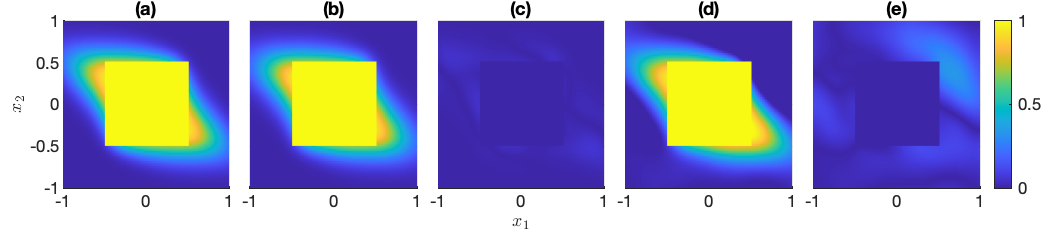}
    \caption{
    (a) Safety probabilities $V_{0}^{\rm DP}(x_{0})$ for a 2-D stochastic chain of integrators computed using dynamic programming over a time horizon $N = 5$. 
    (b) Safety probabilities $V_{0}^{KM}(x_{0})$ computed using \texttt{KernelEmbeddings}
    (c) Absolute error $\vert V_{0}^{\rm DP}(x_{0}) - V_{0}^{KM}(x_{0}) \vert$.
    (d) Safety probabilities $V_{0}^{RFF}(x_{0})$ computed using \texttt{KernelEmbeddingsRFF}
    (e) Absolute error $\vert V_{0}^{\rm DP}(x_{0}) - V_{0}^{RFF}(x_{0}) \vert$.
    }
    \label{fig: double integrator}
\end{figure*}

We first consider a toy example to demonstrate the use of the algorithms and compare the output against known results. 
Consider a $n$-dimensional stochastic chain of integrators \cite{Vinod_Oishi_2020, thorpe2019stochastic}, in which the input appears at the 
$n^{\mathrm{th}}$ 
derivative
and each element of the state vector is the discretized integral of the element that follows it. The
dynamics with sampling time $T$ are given by:
\begin{align}
  x_{k+1} =
  \begin{bmatrix}
    1 & T & \cdots & \frac{T^{n-1}}{(n-1)!} \\
    0 & 1 & \cdots & \frac{T^{n-2}}{(n-2)!} \\
    \vdots & \vdots & \ddots & \vdots \\
    0 & 0 & 0 & 1
  \end{bmatrix}
  x_{k} +
  \begin{bmatrix}
    \frac{T^{n}}{n!} \\
    \frac{T^{n-1}}{(n-1)!} \\
    \vdots \\
    T
  \end{bmatrix}
  u_{k} +
  w_{k}
\end{align}
For the purpose of comparison against other algorithms, we restrict ourselves to the $2$-dimensional case and Gaussian disturbances.
%
We compare the computation time of the algorithms with a single evaluation point $x_{0} = \boldsymbol{0}$ against several other algorithms present in \texttt{SReachTools}. 
The results are displayed in Table \ref{table: computation time}.
Note that the kernel-based algorithms do not compute the \emph{optimal} value functions unless the policy is the \emph{maximally safe} Markov policy $\pi^{*}$ \cite{Summers_Lygeros_2010}. 
This makes direct comparison of the algorithms in \texttt{SReachTools} difficult, since the existing algorithms also perform controller synthesis to obtain the maximal reach probability. 
We then seek to quantify the effect of the sample size $M$ on the computation time. We vary $M$ and plot the computation time as a function of $M$ in Figure \ref{fig: computation time}. This shows that as we increase the sample size $M$, the computation time increases exponentially.

In order to validate the approach, we compute the safety probabilities using dynamic programming $V_{0}^{\rm DP}(x_{0})$. 
We then generate observations of the system by choosing $x_{i} \in \mathcal{X}$, $i = 1, \ldots, 2500$, uniformly in the range $[-1.1, 1.1]^{2}$ and collect a sample $\mathcal{S}$ of size $M = 2500$. 
The dynamic programming results are shown in Figure \ref{fig: double integrator}(a).
Figure \ref{fig: double integrator}(b) shows the safety probabilities $V_{0}^{KM}(x_{0})$ computed for a 2-dimensional integrator system using \texttt{KernelEmbeddings}.
We then treat the dynamic programming solution as a truth model, and plot the absolute error $\vert V_{0}^{\rm DP}(x_{0}) - V_{0}^{KM}(x_{0}) \vert$ in Figure \ref{fig: double integrator}(c).
In order to evaluate \texttt{KernelEmbeddingsRFF}, we then collect a sample of frequency realizations of size $D = 15000$ and compute the safety probabilities $V_{0}^{RFF}(x_{0})$.
The results are shown in Figure \ref{fig: double integrator}(d), where Figure \ref{fig: double integrator}(e) shows the absolute error, $\vert V_{0}^{\rm DP}(x_{0}) - V_{0}^{RFF}(x_{0}) \vert$.

\subsection{Repeated Planar Quadrotor}

\begin{figure}
    \centering
    \includegraphics[width=\columnwidth,height=100px]{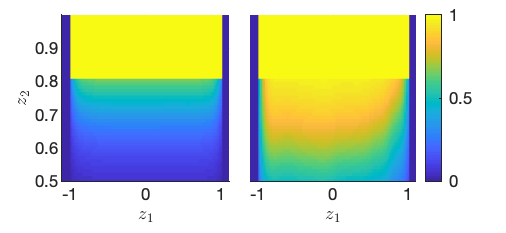}
    \caption{(left) First-hitting time safety probabilities for a single planar quadrotor system with a Gaussian disturbance and (right) with a beta distribution disturbance over the horizon $N = 5$.}
    \label{fig: planar quadrotor}
\end{figure}

This example is used to showcase the ability of the system to handle high-dimensional systems \cite{thorpe2019stochastic}. 
This problem can be interpreted as a simplification of formation control for a large swarm of quadrotors, where we compute the safety probabilities for the entire swarm as they are controlled to reach a particular elevation.
The nonlinear dynamics of a single quadrotor are given by
\begin{align}
    \begin{split}
        m\ddot{x} &= -(u_{1} + u_{2}) \sin(\theta) \\
        m\ddot{y} &= (u_{1} + u_{2}) \cos(\theta) - mg \\
        \mathcal{I}\ddot{\theta} &= r(u_{1} - u_{2})
    \end{split}
\end{align}
where $x$ is the lateral position, $y$ is the vertical position, $\theta$ is the pitch, and we have the constants
intertia $\mathcal{I} = 2$, length $r = 2$,
mass $m = 5$, and
$g = 9.8$ is the gravitational constant. 

The safety probabilities for a single planar quadrotor with a Gaussian and non-Gaussian disturbance are shown in Figure \ref{fig: planar quadrotor}. 
We then formulate the dynamics for a swarm of quadrotors by repeating the dynamics until we have over a million state variables. 
We then generate a sample $\mathcal{S}$ of size $M = 1000$ for the repeated system and compute the safety probabilities using \texttt{KernelEmbeddings}. Then, we collect a sample of frequency realizations of size $D = 15000$ and compute the safety probabilities using \texttt{KernelEmbeddingsRFF} over a single time step $N = 1$. 
The mean computation time for \texttt{KernelEmbeddings} was 1.23 hours, and for \texttt{KernelEmbeddingsRFF} the mean computation time was 44.63 seconds.
This shows that \texttt{KernelEmbeddingsRFF} can be used to compute the safety probabilities for extremely high-dimensional systems.

\subsection{Cart-Pole}

\begin{figure*}
    \centering
    \includegraphics[width=\linewidth,height=100pt]{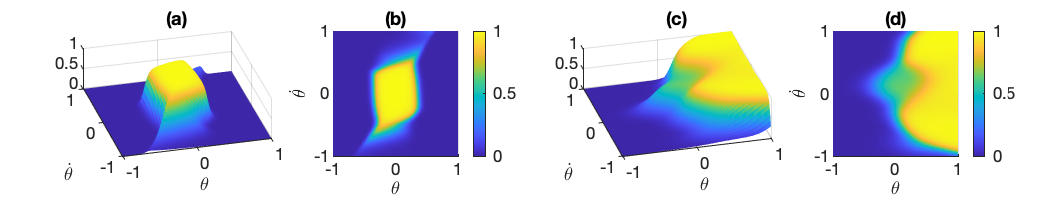}
    \caption{
    (a), (b) 3-D and 2-D cross-sections respectively, of the safety probabilities computed using \texttt{KernelEmbeddings} for linearized Cart-Pole system, and (c), (d) 3-D and 2-D cross-sections respectively, of safety probabilities computed using \texttt{KernelEmbeddings} for the nonlinear cart-pole system.
    }
    \label{fig: cart pole}
\end{figure*}


The following examples are used to showcase the ability of the algorithms to compute the safety probabilities for systems with neural network controllers.

\subsubsection{Linearized Cart-Pole}

The dynamics for the linearized cart-pole system \cite{manzanas2019arch} are given by:
\begin{align}
\begin{split}
\label{eqn: linearized cart pole dynamics}
	\ddot{x} &= 0.0043 \dot{\theta} - 2.75 \theta + 1.94 u - 10.95 \dot{x} \\
	\ddot{\theta} &= 28.58 \theta - 0.044 \dot{\theta} - 4.44 u + 24.92 \dot{x}
\end{split}
\end{align}
We add an additional Gaussian disturbance $\mathcal{N}(0, \Sigma)$ with $\Sigma = 0.01 I$ to the dynamical state equations, which can simulate dynamical uncertainty or minor system perturbations.
The control input is computed via a feedforward neural network controller \cite{manzanas2019arch}, which takes the current state and outputs a real number $u \in \mathbb{R}$, which can be interpreted as the input torque. 

Figure \ref{fig: cart pole}(a) shows a cross-section of the safety probabilities for the system computed using \texttt{KernelEmbeddings} holding $x$ and $\dot{x}$ constant at $0$, and Figure \ref{fig: cart pole}(b) shows a 2-D projection. The algorithms are agnostic to the structure of the dynamics, which means we do not require any prior knowledge of the structure of the neural network controller. Because of this, the algorithms are able to perform verification of systems that incorporate learning enabled components.


\subsubsection{Nonlinear Cart-Pole}

We then analyzed a nonlinear cart-pole system with a neural network controller \cite{manzanas2019arch}, with dynamics given by:
\begin{align}
\begin{split}
\label{eqn: nonlinear cartpole system}
	\ddot{x} ={} & \frac{u + ml\omega^{2} \sin(\theta)}{m_{t}} \\
	&- \frac{ml (g \sin(\theta) - \cos(\theta)) (\frac{u + ml\omega^{2} \sin(\theta)}{m_{t}})}{l(\frac{4}{3} - m \frac{\cos^{2}(\theta)}{m_{t}})} \frac{\cos(\theta)}{m_{t}} \\
	\ddot{\theta} ={} & \frac{g \sin(\theta) - \cos(\theta) (\frac{u + ml\omega^{2} \sin(\theta)}{m_{t}})}{l(\frac{4}{3} - m \frac{\cos^{2}(\theta)}{m_{t}})}
	\frac{\cos(\theta)}{m_{t}}
\end{split}
\end{align}
where $g = 9.8$ is the gravitational constant, the pole mass is $m = 0.1$, half the pole's length is $l = 0.5$, and $m_{t} = 1.1$ is the total mass.
The control input, $u \in \lbrace -10, 10 \rbrace$, which affects the lateral position of the cart, is chosen by the neural network controller \cite{manzanas2019arch}. 
We add an additional Gaussian disturbance $\mathcal{N}(0, \Sigma)$ with $\Sigma = 0.01 I$ to the dynamical state equations.

Figure \ref{fig: cart pole}(c) shows a 3-D representation of the safety probabilities for the system computed using \texttt{KernelEmbeddings}, and Figure \ref{fig: cart pole}(d) shows a 2-D projection. 
This example shows that we can handle systems which are traditionally very difficult to model and analyze, such as nonlinear systems and systems with neural network controllers.


\section{Conclusion \& Future Work}
\label{section: conclusion}

We presented a data-driven stochastic reachability module for \texttt{SReachTools} based on conditional distribution embeddings. 
These algorithms add to the suite of stochastic reachability algorithms already present in the toolbox and enable point-based stochastic reachability for a wide variety of stochastic systems.
We would like to extend the kernel-based algorithms to enable model-free controller synthesis and broaden the capability of the algorithms to handle a wider variety of systems, such as hybrid system models. 


\bibliographystyle{IEEEtran}
\bibliography{IEEEabrv, shortIEEE, bibliography}

\end{document}